\documentclass[10]{article}

\def\ba{\mathbf{a}}   
\def\bb{\mathbf{b}}   
\def\be{\mathbf{e}}   
\def\bff{\mathbf{f}}   
\def\bi{\mathbf{i}}   
\def\bj{\mathbf{j}}   
\def\bk{\mathbf{k}}   
\def\bv{\mathbf{v}}   
\def\bw{\mathbf{w}}   
\def\bN{\mathbf{N}}   
\def\bW{\mathbf{W}}   

\def\C{\mathbb{C}} 
\def\G{\mathbb{G}} 
\def\H{\mathbb{H}}  
\def\Q{\mathbb{Q}}
\def\R{\mathbb{R}}  

\def\no{\noindent}

\def\beq{\begin{equation}}
\def\eeq{\end{equation}}

\def\w{\wedge}

\usepackage{array}
\usepackage{tabularx}
\usepackage{amsfonts}
\usepackage{amssymb}
\usepackage{tikz-cd}
\usetikzlibrary{matrix,arrows,decorations.pathmorphing}

\def\bpm{\begin{pmatrix}}
\def\epm{\end{pmatrix}}
\makeindex           
\newtheorem{thm1}{Theorem}
\newtheorem{thm2}[thm1]{Theorem}
\newtheorem{thm3}[thm1]{Theorem}
\title{Periodic Table of Geometric Numbers}
\author{Garret Sobczyk \\
Universidad de las Am\'ericas-Puebla \\
 Departamento de Actuaría F\'isica y Matem\'aticas \\
72820 Puebla, Pue., M\'exico}
	
\begin{document}

\maketitle

\begin{abstract} 
 Perhaps the most significant, if not the most important, achievements in chemistry and physics are the Periodic Table of the Elements in Chemistry and the Standard Model of Elementary Particles in Physics. A comparable achievement in mathematics is the Periodic Table of Geometric Numbers discussed here. In 1878 William Kingdon Clifford discovered the defining rules for what he called geometric algebras. We show how these algebras, and their coordinate isomorphic geometric matrix algebras, fall into a natural periodic table, sidelining the superfluous definitions based upon tensor algebras and quadratic forms. 
 
\smallskip
\no {\em AMS Subject Classification:} 15A63, 15A66, 81R05, 01A55
\smallskip

\no {\em Keywords: Clifford algebra, geometric matrices, quaternions, }
\end{abstract}

 \section*{0\quad Introduction}
 
 The periodic table of elements, first recognized as such, was developed by 
 Dmitri Mendeleev and published in 1869. Since that time it has undergone many extensions and has proved itself to be of fundamental importance in research and predicting chemical reactions. The most standard version has 7 rows, called periods, and 18 columns, called groups. The period in which an element belongs depends upon the structure of its electrons in their orbital shells, and all the elements in a particular group share similar properties \cite{hptable}.
 
 The standard model of particle physics consists of two classes of particles depending upon their spin. Fermions consisting of 6 quarks and 6 leptons all have spin 1/2, and force bosons consisting of 5 particles, including the familiar photon, have spin 0. The standard model has enjoyed huge successes in making experimental predictions, and models quantum field theory, \cite{ahey}.
 
 It would be a mistake to talk about the development of the Periodic Table of Chemistry and the Standard Model of Particle Physics without mentioning the conceptual developments that were simultaneously taking place in mathematics, developments that often made possible the revolutionary progress in chemistry, physics and other areas \cite{unreason}. Particularly important were the developments of the real and complex number systems, and their higher dimensional generalizations to Hamilton's quaternions (1843), Grassmann algebras (1844), and Clifford's geometric algebras (1878), \cite{wkc,WKC1882}, \cite{grass1862}, \cite{wrh}.
 Unfortunately, this early flowering of mathematics was made before the appearance of Einstein's four dimensional Theory of  Special Relativity in 1905. 
 
 As a consequence of this fluke of history, elementary mathematics has been hobbled by the 3-dimensional Gibbs Heaviside vector analysis developed just prior to Einstein's famous theory, to the neglect in the development of Grassmann's and Clifford's higher dimensional geometric number systems. The main purpose of this article is to show how Clifford's geometric algebras, as demonstrated by the Periodic Table of Geometric Numbers, is the culmination of the millennium old development of the concept of number itself. It is a basic structure in linear algebra, and need not be shackled to  unnecessary concepts from higher mathematics. The prerequisites for understanding of this article is a knowledge of matrix multiplication, and the concept of a vector space \cite{hmatrix}.
 
 The importance of geometric algebras in physics was first recognized by 
 Brauer and Weyl \cite{b/w1935}, and Cartan \cite{cartan1966}, particularly in connection with the concept of $2$- and $4$-component spinors at the heart of the newly minted {\it quantum mechanics} \cite{sources1967}. The concept of a spinor arises naturally in Clifford algebras, and  much work has been carried out by mathematicians and physicists since that time \cite{dss1992}. 
  More recently, the importance of Clifford geometric algebras has been recognized in the computer science and engineering communities, as well as in efforts to develop the mathematics of quantum computers \cite{HD02}. 
  Over the last half century Clifford's geometric algebras have become a powerful geometric language for the study of the atomic structure of matter in elementary particle physics, Einstein's special and general theories of relativity, string theories and super symmetry. In addition, geometric algebras have found their way into many areas of mathematics, computer science, engineering and robotics, and even the construction of computer games, \cite{gs2020,Sob2012,S0}. 
 
 \section{What is a geometric number?}
 
 The development of the real and complex number systems have a long and torturous history, involving many civilizations, wrong turns, and dead ends. It is only recently that any real perspective on the historical process has become possible. The rapid development of the geometric concept of number, over the last 50 years or so, is based upon Clifford's seminal discovery of the rules of geometric algebra in 1878. 
 Those rules are succinctly set down in the following
 
 \begin{quote}{\bf Axiom:} {\em   The real number system can be geometrically extended to include new, anti-commutative square roots of $\pm 1$, each new such root representing
 the direction of a unit vector along orthogonal coordinate axes of a Euclidean or
 pseudo-Euclidean space $ \R^{p,q}$, where $p$ and $q$ are the number of new square roots of $+1$ and $-1$, respectively.} 
 \end{quote}
 
    The resulting real geometric algebra, denoted by
 \[  \G_{p,q} := \R(\be_1, \ldots, \be_p, \bff_1, \ldots, \bff_q), \]
 has dimension $2^{p+q}$ over the real numbers $\R$, and is said to be {\it universal} since no further relations between the new square roots are assumed.
 Since the elements in $\G_{p,q}$ satisfy the same rules of addition and multiplication as the addition and multiplication of matrices, it is natural
 to consider matrices with entries in $\G_{p,q}$. Indeed, as we will shortly see, the elements of a geometric algebra provide a natural geometric basis for matrix algebra,  and, taken together, form an integrated framework which is more
 powerful than either when considered separately, \cite{SNF,S08}.
 Also considered are complex geometric algebras and their corresponding complex matrix algebras, \cite[p.75]{matrixG}. 
 
 A geometric number $g\in \G_{p,q}$ can be decomposed into its various $k$-vector parts. Letting $n=p+q$, 
    \beq g=\sum_{k=0}^{n} \langle g \rangle_k, \label{k-part} \eeq
  where the $k$-vector part 
  \[  \langle g \rangle_k:=\sum_{i=1}^{\pmatrix{n \cr k}}g_{i_1 \cdots ,i_k}\be_{i_1,\cdots,i_r}\bff_{i_{r+1},\cdots i_k}\in \G_{p,q}^k  ,\]
  for $1 \le i_1 < \cdots <i_r \le p$ and $1 \le i_{r+1} < \cdots <i_k \le q$, and
    \[ \be_{i_1,\cdots,i_r}:=\be_{i_1}\cdots \be_{i_r}\ \ {\rm and} \ \ \bff_{i_{r+1},\cdots,i_k}:=\bff_{i_{r+1}}\cdots \bff_{i_k}.        \]
  The coordinates $g_{i_1, \cdots ,i_k}$ are either real or complex numbers, depending upon whether $\G_{p,q}$ is a {\it real} or {\it complex} geometric algebra. 
 The {\it standard basis} of the real geometric algebra,
 \beq \G_{p,q}:=span_\R \{1,\be_1,\ldots, \be_p,\bff_1, \ldots, \bff_q, \ldots, \be_{1,\ldots,p}\bff_{1,\ldots,q}   \},   \label{stdbasis} \eeq
 has dimension
 $ 2^n= \sum_{k=0}^n \pmatrix{n\cr k}$.
 The element of highest grade, the $n$-vector $\be_{1,\ldots,p}\bff_{1,\ldots,q}\in \G_{p,q}^n $ is called the {\it pseudoscalar} of $\G_{p,q}$.
 
  There are 3 basic {\it conjugations} defined on $\G_{p,q}$, {\it reverse, inversion} and {\it mixed}. For any $g\in \G_{p,q}$, the {\it reverse}  $g^\dagger$ of $g$ is obtained by reversing the order of the products of vectors which define $g$. For example, if 
 \[ g=1+2\be_1-3\bff_{12}+4\be_{12}\bff_3 , \ \ {\rm then} \ \ g^\dagger = 1+2\be_1+3\bff_{12}-4\be_{12}\bff_3 .\]
  The
 {\it inversion} $g^-$ of $g$ is obtained by replacing every vector in $g$ by its negative. So for   
 \[ g=1+2\be_1-3\bff_{12}+4\be_{12}\bff_3 , \ \ {\rm then} \ \ g^- = 1-2\be_1+3\bff_{12}-4\be_{12}\bff_3 .\]
 Finally, the {\it mixed conjugation} $g^*:={g^\dagger}^-$ of $g$ is the composition of the reversion and inversion of $g$. For
 \[ g=1+2\be_1-3\bff_{12}+4\be_{12}\bff_3 , \ \ {\rm then} \ \ g^* :={g^\dagger}^-= 1-2\be_1+3\bff_{12}+4\be_{12}\bff_3 .\]
  
 This brief summary of the definition and defining properties of a geometric algebra is not complete without mentioning the definitions of the {\it inner and outer products} of $k$-vectors, and the {\it magnitude} of a geometric number. The most famous formula in geometric algebra is the decomposition of the geometric product of two vectors into a symmetric {\it inner product} and a skew-symmetric {\it outer product}. For $\bv, \bw \in \G_{p,q}^1$, 
   \[ \bv \bw = \frac{1}{2}( \bv \bw + \bw \bv)+ \frac{1}{2}( \bv \bw - \bw \bv)=\bv\cdot \bw + \bv \w \bw ,  \]
 where the scalar inner product $\bv \cdot \bw:= \frac{1}{2}( \bv \bw + \bw \bv)=\langle \bv \bw \rangle_0$, and bivector-valued outer product
    $\bv \w \bw:= \frac{1}{2}( \bv \bw - \bw \bv)=\langle \bv \bw \rangle_2$. 
    
    More generally, for a vector $\bv$ and a $k$-vector $\bW_k$,
    \[ \bv \bW_k = \frac{1}{2}\big( \bv \bW_k +(-1)^{k+1} \bW_k \bv\big)+ \frac{1}{2}\big( \bv \bW_k - (-1)^{k+1} \bW_k \bv\big) ,  \] 
where 
\[ \bv \cdot \bW_k:= \frac{1}{2}\big( \bv \bW_k +(-1)^{k+1} \bW_k \bv\big) = \langle \bv \bW_k \rangle_{k-1}\] 
and 
\[ \bv \w \bw:= \frac{1}{2}\big( \bv \bW_k - (-1)^{k+1} \bW_k \bv\big)= \langle \bv \bW_k \rangle_{k+1}. \]
 A general {\it inner product} between geometric numbers $g_1,g_2 \in \G_{p,q}$ is defined by 
 \[  g_1 * g_2 :=\langle g_1 g_2^\dagger \rangle_0,\]
  and the {\it magnitude} of $g\in \G_{p,q}$ is 
   \[ |g|:=|\langle g g^\dagger\rangle_0|^{\frac{1}{2}}. \] 
 A complete treatment of geometric algebra, and the myriad of algebraic identities, can be found in the books \cite{H/S,SNF}.

 The periodic table of real geometric algebras, shown in Table \ref{tab:table1}, consists of 8 rows and 15 columns. Each of the 36 geometric algebras in the Table is specified by an isomorphic {\it geometric matrix algebra} of coordinate matrices over one of the five fundamental building blocks. By giving a direct construction of each geometric algebra from its  isomorphic coordinate matrix algebra, we have a tool for carrying out and checking calculations in terms of the well-known ubiquitous matrix addition and multiplication. It is amazing that geometric algebras today, 142 years after their discovery, are still not recognized as the powerful higher dimensional geometrical number systems that they are, with applications across the scientific, engineering and computer science and robotics communities. 
 
  \begin{table}
  	\centering
  	\caption{Classification of real geometric algebras $\G_{p,q}$.}
  	\label{tab:table1}
  	
  	{\footnotesize 
  		\[ \begin{tabular*}{\columnwidth}[t]{@{}c@{}c@{}c@{}c@{}c@{}c@{}c@{}c@{}c@{}c@{}c@{}c@{}@{}c@{}c@{}c@{}}
  		7 & 6  & 5&4 &3&2&1&0 &-1&-2& -3&-4&  -5&-6 &-7 \\
  		0$\quad \quad\quad \quad$    && & &  & & & $\R$ & & & & & &&  \\
  		1$\quad \quad\quad \quad$ &&&&&& $^2{\R}$ &&  ${\C}$ &&&&&& \\
  		2$\quad \quad\quad \quad$  &&&&& $M_2\R$ && $M_2$$\R$ && 
  		$\Q$ & & &&& \\
  		3$\quad \quad\quad \quad$  &&&& $M_2$${\C}$ && $M_2$$^2{\R}$ &&  $M_2{\C}$ && $^2\Q$ &&&& \\
  		4$\quad \quad\quad \quad$ &&& $M_2$$\Q$&& $M_4$$\R$&&$M_4$$\R$&&$M_2$$\Q$&& $M_2$$\Q$\\
  		5$\quad \quad\quad \quad$ &&$ M_2$$^2\Q$ && $M_4{\C}$&&$M_4$$^2{\R}$ && $M_4{\C}$&& $M_2$$^2\Q$&& $M_4{\C}$ \\
  		6$\quad \quad\quad \quad$ & $M_4$$\Q$ && $M_4$$\Q$ && $M_8$$\R$&&$M_8$$\R$ && $M_4$$\Q$&& $M_4$$\Q$&& $M_8$$\R$ \\
  		7$\quad \ $  $M_8\C$ && $M_4$$^2\Q$ && $M_8\C$ && $M_8$$^2{\R}$&&$M_8{\C}$ && $M_4$$^2\Q$&& $M_8{\C}$&& $M_8$$^2{\R}$ \\
  		$i^2 \quad \quad \  -$ & $-$ &$+$ &$ +$ & $-$ & $-$ & $+$ & $+$ &$ -$ & $-$ &$ +$ &$ +$ &$ -$ &$ -$ &$ +$ \\
  		\end{tabular*}  \] }
  \end{table}
 
  \begin{table}[h!]
 	\centering
 	\caption{Budinich/Trautman Clifford Clock.}
 	\label{tab:table2} 
 	\bigskip
 	{\footnotesize
 		\begin{tikzcd}
 			&& \R  \arrow{rd}
 			\\ & \,^2{\R}  \arrow{ru}\arrow[dashed]{rd}   & &  \arrow{rd} \C \\
 			{\R}  \arrow{ru}  &&\arrow[dashed]{rd} && \Q \arrow{ld}
 			\\ &  {\C}   \arrow{lu}   & &   ^2\Q \arrow{ld}  \\
 			&& \Q \arrow{lu}  
 		\end{tikzcd} }
 	\end{table}
 	
 	The Clifford clock, Table \ref{tab:table2}, contains in coded form exactly the same information about the geometric matrix representation of $\G_{p,q}$ as in Table \ref{tab:table1}, \cite{budtraut}. Starting from the top $\R$ of each table, to get to $\G_{p,q}$ in Table \ref{tab:table1}, one takes $0 \le q \le 7$ successive steps to the right down, followed by $0\le p \le 7-q$ steps to the left down. The corresponding entry on the Clifford clock in Table \ref{tab:table1}, is found by advancing $q$ clockwise Clifford hours followed by $p$ counterclockwise Clifford hours. The total number of hours elapsed, $n=p+q$, determines the dimension of the matrix algebra over the block hour entry on the Clifford Clock reached. For the geometric algebra $\G_{3,2}$ in Table \ref{tab:table1}, one proceeds from the top $\R$ two steps to the right and down, followed by three steps to left and down to reach the entry $M_4(\,^2\R)$. On the Clifford clock, we proceed clockwise two hours to reach the block hour $\Q$, and then three more hours counterclockwise to reach the block hour $\,^2\R$. The number $n=2+3=5$ $\R$-block hours corresponds the $4\times 4$ real matrix algebra $M_4(\,^2\R)$ isomorphic to the $2\cdot 2^2\times 2^2=32$ dimensional real geometric algebra $\G_{3,2}$.

\section{Building blocks}

\begin{figure}
  	\includegraphics[width=0.8\linewidth]{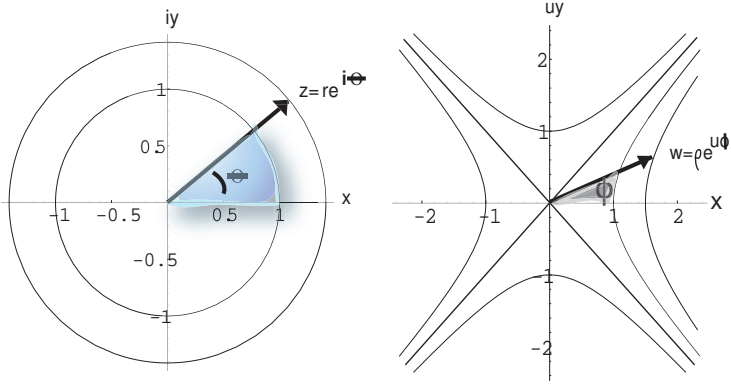}
  	  			\caption{a) In the $r$-circle the shaded area $=\frac{\theta}{2}$. b) In the $\rho$-hyperbola the shaded area $=\frac{\phi}{2}$. }
 		\label{rcircle}
 \end{figure}

    The two most familiar building blocks are the real and complex number systems, denoted by $\R$ and $\C$, respectively. Geometrically, the real numbers $\R\equiv \G_{0,0}$ are pictured on a horizontal line, not shown, with the 0 point at its center. The complex numbers $\C:=\G_{0,1}=\R(i)$, in turn, are pictured in the complex number plane, shown in Figure \ref{rcircle}: a), with the real $x$-axis, and the {\it imaginary} $iy$-axis. The hyperbolic number plane $\H:=\G_{1,0}= \R(u) $, the twin sister of the complex number plane, is pictured in Figure \ref{rcircle}: b). It has a real $x$-axis and a {\it unipotent} $uy$-axis. Whereas $i^2=-1$ for the imaginary unit $i$, the less familiar hyperbolic unit $u\notin \R$ has square $u^2=+1$, \cite[Ch.1]{matrixG}.

    William Hamilton (1788-1856) invented the quaternions $\Q$, consisting of three
    new {\it anticommutative} numbers $\{\bi, \bj, \bk  \}$, and satifying the rules
    \beq   \bi^2 = \bj^2=\bk^2=-1,  \ {\rm and} \ \ \bi \bj  \bk = -1. \label{quat} \eeq  
    Because the quaternions $\bi \bj \bk=-1$, they do not satisfy the {\it universal property} of a geometric algebra; hence $\R(\bi,\bj, \bk)$ is not a geometric algebra. 
    
    Hamilton interpreted his unit quaternions to be vectors along the $x,y,z$-axes.
    However, Hamilton's rules (\ref{quat}) for his anticommutative units coincide exactly with those of the geometric algebra 
    \[ \G_{0,2}:=\R(\bff_1,\bff_2) =span_\R \{ 1,\bff_1,\bff_2,\bff_1 \bff_2\},  \]
    with $\bi =\bff_1, \bj=\bff_2,\bk = \bff_{1}\bff_2$. The geometric algebra 
    $\G_{0,2}$ satisfies the universal property, since 
    $\bff_1\bff_2 \in \G_{0,2}^2$ is a bivector independent of the vectors $\bff_1,\bff_2$.  
    
    Whereas Hamilton identified his quaternions with vectors, we see in $\G_{0,2}$ that the anticommutative quaternions can equally well be identified with a mixture of the vectors $\bff_1, \bff_2$ with the bivector $\bff_{12}=\bff_1 \bff_2$. The important point is that the anticommutative relations themselves are often more important than the particular geometric interpretation of being a vector or a bivector, or something else. Thus, we can
    justifiably make the definition $\Q:=\G_{0,2}$.   
    
    So far we have identified $\R , \C=\G_{0,1}, \H = \G_{1,0} $ and $Q=\G_{0,2}$, four of the five building blocks of real geometric algebras. We have already noted that extending the real numbers $\R$ by the unipotent $u\notin \R$ gives the geometric algebra $\G_{1,0}$. Indeed, extending any geometric algebra $\G_{p,q}$ by a {\it new number}
    $u$, which commutes with all of the numbers in $\G_{p,q}$, has the effect of
    doubling the geometric algebra into the {\it double geometric algebra}, denoted by
    \beq ^2\G_{p,q}:=\G_{p,q}(u)=\{ (g_1,g_2)| \ {\rm for \ all \ g_1,g_2 \in \G_{p,q}}\}, \label{doubleG}  \eeq
    where addition and multiplication in $\G_{p,q}(u)$ is naturally defined by
    \[ (g_1,g_2)+(g_3,g_4)=(g_1+g_3,g_2+g_4), \ \ {\rm and } \ \  (g_1,g_2)(g_3,g_4)=(g_1g_3,g_2g_4), \]
    respectively. 
    
      Note that the hyperbolic numbers are just double real numbers, i.e.,
      $\H:=\, ^2\R = \, ^2\G_{1,0}$. To see this, we define 
      \[u_+:=\frac{1}{2}(1 + u), \ {\rm and} \ u_-:=\frac{1}{2}(1-u), \]
 and verify that $u_\pm^2=u_\pm^2, u_+u_-=0$, so that $u_+$ and $u_-$ are {\it mutually annihiliating idempotents}, and $u_+-u_-=u$. It then follows that any number $X=x_1+x_2 u \in \R(u)$ for $x_1,x_2 \in \R $, can be written in the form
 \[ X =X(u_++u_-)=Xu_++Xu_-=(x_1+x_2)u_++(x_1-x_2)u_-, \]
 or, equivalently, $X=(x_1+x_2,x_1-x_2)\in \, ^2\R$. 
 
 With the concept of a {\it double geometric algebra} (\ref{doubleG}), we can identify 
 all six building blocks of real and complex geometric algebras,
 \beq  \R,\C=\R(\bff_1) , \Q=\R(\bff_1,\bff_2),\,^2\Q=\R(\bff_1,\bff_2,\bff_3) , \,^2\R=\R(u),  \,^2 C=\C(u) . \label{buildingB} \eeq 
  To see that 
  \beq \,^2\Q=\G_{0,3}=\R(\bff_1,\bff_2,\bff_3) =\R(\bff_1,\bff_2)(\bff_3)= \G_{0,2}(\bff_3) , \label{efnotation} \eeq 
   note that $\bff_{123}$ is in the center of $Z(\G_{0,3} )$ since it commutes with $\bff_1,\bff_2$ and $\bff_3$. Since, in addition, $\bff_{123}^2=+1$, it can play the role of $u$ so that
  $\,^2\Q = \Q(\bff_{123})$.   
  Examining Table \ref{tab:table1}, we see that all real geometric algebras can be represented in terms of isomorphic matrix algebras over five of the six building blocks. Table \ref{tab:table2} shows how these five building blocks define the
  Clifford Clock, which contains the same information as Table \ref{tab:table1} but in coded form. The sixth building block, $\,^2\C$, will be used later when constructing the complex geometric algebras $\G_{n}(\C)$. 
  
  In (\ref{efnotation}), we have also introduced the idea of extending a geometric algebra by an additional anticommuting vector which has square $-1$. More generally, given a geometric algebra $\G_{p,q}$, we define the {\it f-extension}
  \beq \G_{p,q+1}=\G_{p,q}(\bff_{q+1}):=\R(\be_1,\ldots,\be_p,\bff_1,\dots,\bff_q ,\bff_{q+1})\label{fextention}    \eeq  
  and the {\it e-extension}  
   \beq \G_{p+1,q}=\G_{p,q}(\be_{p+1}):=\R(\be_1,\ldots,\be_p,\be_{p+1},  \bff_1,\dots,\bff_q).\label{eextention}    \eeq  
    \section{Geometric algebras $\G_{1,1}$ and $\G_{1,2}$} 

\begin{figure}[h]
	\includegraphics[width=0.40\linewidth]{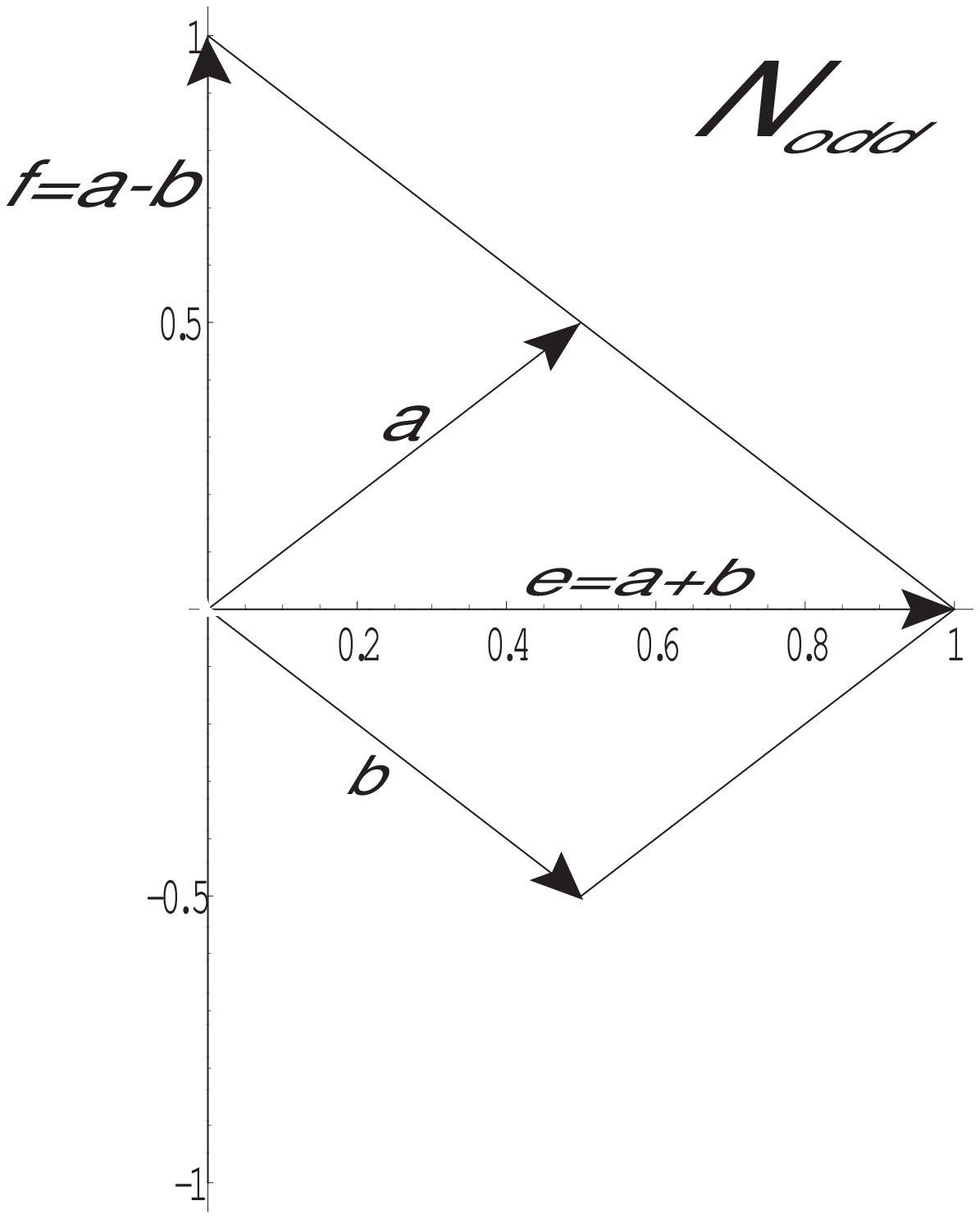}
	\includegraphics[width=0.40\linewidth]{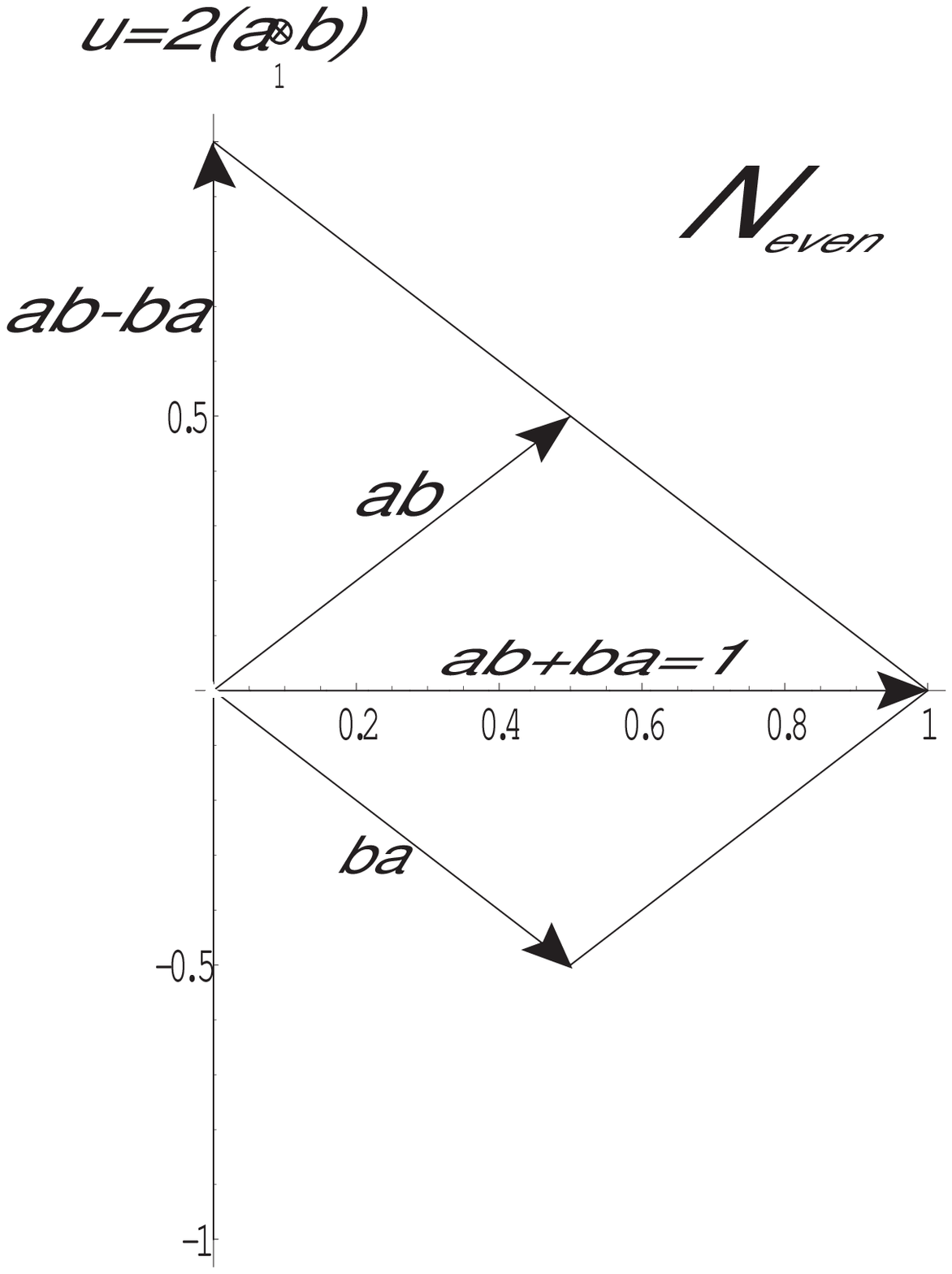}
	\caption{The odd g-number plane $\bN_{odd}$.\quad The even g-number plane $\bN_{even}$.}
	\label{evenoddhyp}
\end{figure}

In the standard basis, the geometric algebra 
\beq \G_{1,1}:=\R(\be_1,\bff_1) = span_\R\{1 , \be_1,\bff_1,\be_1 \bff_2 \}. \label{stdbasisG11} \eeq
The {\it odd} part $\bN_{odd}\subset\G_{1,1}$ consists of the vectors of $\G_{1,1}$, and the {\it even} part $\bN_{even}:=\G_{1,1}^+$, consists of the real numbers $\R$ and bivectors of $\G_{1,1}$.  
The standard {\it null vector basis} of the geometric algebra $\G_{1,1}$ is defined in terms of the {\it null vectors} $\{\ba , \bb \}$.
We have
  \beq   \G_{1,1}:=\R(\ba, \bb):=span_\R\{\ba,\bb,\ba \bb, \bb \ba \} , \label{nullvecbasis} \eeq
  where

\beq \ba := \frac{1}{2}(\be_1+\bff_1), 
\ {\rm and} \ \bb :=\frac{1}{2}(\be_1 -\bff_1). \label{nullabdef} \eeq
The even and odd parts of $\G_{1,1}$, together with the basis vectors are pictured in Figure \ref{evenoddhyp}.

The properties of the null vectors $\ba$ and $\bb$ are fully characterized by the simple rules
  \begin{itemize}
  	\item[N1)] $\ba^2 = 0 = \bb^2$. {\it The vectors $\ba $ and $\bb $ are null vectors.}
  	\item[N2)] $ \ba \bb + \bb \ba =1 \ \  \iff \ \ \ba \cdot \bb =\frac{1}{2}(\ba \bb + \bb \ba )= \frac{1}{2}$.
  \end{itemize} 
 A Multiplication Table for these null vectors, given in Table \ref{table3}, can be checked by the reader, \cite[p.74]{matrixG}.
  
 \begin{table}[h]
 	\begin{center}
 		\caption{Multiplication Table}
 		\label{table3}
 		\begin{tabular}{c|c|c|c|c} 
 			& $\ba $  & $ \bb $  & $\ba \bb $  & $ \bb \ba $ \cr 
 			\hline
 			$ \ba$  &   0   & $ \ba \bb $     & $  0 $         & $ \ba $ \cr 
 			\hline         
 			$ \bb $  &  $ \bb \ba $  & $ 0 $     & $  \bb $         & $ 0 $ \cr 
 			\hline     
 			$ \ba \bb $    &  $\ba $    &  0  &  $\ba \bb $    & 0  \cr
 			\hline
 			$ \bb \ba $  &  0   &  $\bb $   &  0     & $\bb \ba$   \cr
 		\end{tabular}
 	\end{center}
 \end{table}  The null vector basis (\ref{nullvecbasis}) is particularly powerful because it ties a geometric number $g\in \G_{1,1}$ directly to its {\it geometric matrix} $[g]:=\pmatrix{g_{11} & g_{12} \cr g_{21}& g_{22}}$ of coordinates $g_{ij}\in \R$. Rewriting the null vector basis (\ref{nullvecbasis}) of $\G_{1,1}$ in the matrix form
 \beq  \pmatrix{\bb \ba \cr \ba} \pmatrix{\bb \ba & \bb}= \pmatrix{\bb \ba & \bb \cr  \ba & \ba \bb}, \label{nullstdbasis} \eeq
 we have 
 \beq g =  \pmatrix{\bb \ba & \ba} [g] \pmatrix{\bb \ba \cr \bb} = g_{11}\bb \ba + g_{12}\bb + g_{21}\ba + g_{22}\ba \bb.    \label{matrixg11} \eeq

The equation (\ref{matrixg11}) establishes an {\it algebra isomorphism} $\G_{1,1} \widetilde= M_2(\R )$ between the geometric algebra $\G_{1,1}$ and the matrix algebra $M_2(\R)$. To see this suppose that $g_1 =  \pmatrix{\bb \ba & \ba} [g_1] \pmatrix{\bb \ba \cr \bb} $ and $g_2=  \pmatrix{\bb \ba & \ba} [g_2] \pmatrix{\bb \ba \cr \bb}$. Then
\[ g_1 g_2= \pmatrix{\bb \ba & \ba} [g_1] \pmatrix{\bb \ba \cr \bb}  \pmatrix{\bb \ba & \ba} [g_2] \pmatrix{\bb \ba \cr \bb}      \]  
\[   = \pmatrix{\bb \ba & \ba} [g_1] \pmatrix{\bb \ba  & 0 \cr 0 & \bb \ba } [g_2] \pmatrix{\bb \ba \cr \bb}   = \pmatrix{\bb \ba & \ba} [g_1] [g_2] \pmatrix{\bb \ba \cr \bb} ,      \]
so the product $g_1g_2 \in \G_{1,1}$ correponds to the matrix product $[g_1g_2]\in M_2(\R )$.

 There is also a matrix form of the standard basis (\ref{stdbasisG11}), 
 \[  \G_{1,1}:= \pmatrix{\bb \ba \cr \ba} \pmatrix{\bb \ba & \bb}=\pmatrix{1 \cr \ba}\bb \ba \pmatrix{1 & \bb} \]
 \beq =\pmatrix{1 \cr \ba +\bb}\bb \ba \pmatrix{1 &\ba + \bb} = \pmatrix{u_+ & \be_1u_- \cr  \be_1u_+ & u_-}, \label{mstdbasisG11}  \eeq
 where $u_\pm = \frac{1}{2}(1\pm u \bf)$ for $u=\bb \ba - \ba \bb=2\bb \w \ba= \be_1\bff_1$. In this basis $g\in \G_{1,1}$ takes the form
 \[  g = \pmatrix{ 1 & \be_1}u_+[g]\pmatrix{1 \cr \be_1}=g_{11}u_++g_{12}\be_1 u_-+g_{21}\be_1u_++g_{22}u_- \]
 \[   = \frac{1}{2}(g_{11}+g_{22})+ \frac{1}{2}(g_{21}+g_{12})\be_1+ \frac{1}{2}(g_{21}-g_{12})\bff_1+ \frac{1}{2}(g_{11}-g_{22})\be_1\bff_1 . \]     
 
 The {\it real} geometric algebra, $\G_{1,2}:=\G_{1,1}(\bff_2)=\R(\be_1,\bff_1,\bff_2) $,  is spanned by the standard basis
 \[ \G_{1,2}:=span_\R \{1,\be_1,\bff_1,\bff_2,\be_1 \bff_1,\be_1\bff_2,\bff_{12},\be_1\bff_{12}\} . \]
 The pseudoscalar $i:=\be_1 \bff_{12}$ has the special property that it commutes with the basis vectors $\be_1,\bff_1,\bff_2\in \G_{1,2}^1$, and therefore with all of the elements of $\G_{1,2}$. It follows that $i\in {\cal Z}(\G_{1,2} ) $, the center of the algebra. It is easily checked that 
 \[  i^2=\be_1 \bff_{12}\be_1 \bff_{12}=\be_1^2 \bff_{12}^2 = -1.   \]  
    
    Because $i\in {\cal Z}(\G_{1,2} )$, $i^2=-1$, and $\G_{1,1}\widetilde{=}M_2(\R)$, it follows that $\G_{1,2}\widetilde{=}M_2(\C)$. Thus, the real geometric algebra $\G_{1,2}$ is isomorphic to the {\it complex} matrix  algebra $M_2(\C)$. For $g\in \G_{1,2}$, using (\ref{matrixg11}),
     \beq g =  \pmatrix{\bb \ba & \ba} [g]_\C \pmatrix{\bb \ba \cr \bb} = g_{11}\bb \ba + g_{12}\bb + g_{21}\ba + g_{22}\ba \bb,    \label{cmatrixg11} \eeq 
    for the complex matrix $[g]_\C=[g_{ij}]$ for $g_{ij}\in \C$. Note also that $i:=\be_1\bff_{12}$, implies that $\bff_2=i\be_1\bff_1 $. 
    
    \section{Classification theorems of geometric algebras}
    
    In order to fully understand the Classification Tables of geometric algebras, three  structure theorems are required. The First Structure Theorem establishes an isomorphism between a geometric algebra $\G_{p,q}$ and the even subalgebra $\G_{p,q+1}^+$ of $\G_{p,q+1}$. The second part of the Theorem establishes that
    $\G_{p+1,q}\widetilde{=}\G_{q+1,p}$. 
    
    \begin{thm1} a) The geometric algebra $\G_{p,q} \widetilde=\, \G_{p,q+1}^+ $.
    	
    	b) The geometric algebra $\G_{p+1,q}\widetilde=\, \G_{q+1,p}$.  
    	\end{thm1}
   
   \no {\bf Proof:} a) By definition, 
   \[  \G_{p,q+1}:= \R(\be_1,\ldots, \be_p,\bff, \bff_1,\ldots, \bff_q).       \]
   It is also true that
   \[ \G_{p,q+1}\widetilde{=} \R(\bff \be_1,\ldots, \bff \be_p,\bff,\bff \bff_1,\ldots, \bff \bff_q).      \]	
   It follows that
   \[ \G_{p,q}\widetilde{=} \R(\bff \be_1,\ldots, \bff \be_p,\bff \bff_1,\ldots,\bff \bff_q) = \G_{p,q+1}^+ .      \]
   
   b) By definition, 
   \[  \G_{p+1,q}:= \R(\be ,\be_1,\ldots, \be_p, \bff_1,\ldots, \bff_q).       \]
   Also,
   \[ \G_{q+1,p}\widetilde{=} \R(\be,\be \bff_1,\ldots, \be  \bff_q,\be \be_1,\ldots, \be \be_p).      \]	
   It follows that
   $ \G_{q+1,p} \widetilde = \G_{p+1,q} , $
   since
    \[ \R(\be ,\be_1,\ldots, \be_p, \bff_1,\ldots, \bff_q) = \R(\be,\be \bff_1,\ldots, \be  \bff_q,\be \be_1,\ldots, \be \be_p) 	.\]
    
    $\hfill \square $
   
  \ Note in the proof of the Theorem the re-identification, or regrading, of the geometric algebras that are involved. This is crucial in understanding the algebra isomorphisms between geometric algebras expressed in the Classification Tables \ref{tab:table1} and \ref{tab:table2}. For example, the result $\G_{p+1,q}\widetilde= \G_{q+1,p}$, tells us that the geometric algebra $\G_{p+1,q}$ at the vertex $(p+1+q,p+1-q)$ is isomorphic to the geometric algebra $\G_{q+1,p}$ at the vertex $(q+1+p,q+1-p)$ in Table \ref{tab:table1}, respectively. In particular, for $p=4,q=0$, the result shows that $\G_{5,0}\widetilde{=}\G_{1,4}$.   
    
  The Second Structure Theorem relates the higher dimensional geometric algebra $\G_{p+k,q+k}$  to a $2^k\times 2^k$ matrix algebra over the geometric algebra $\G_{p,q}$. It allows us to move vertically between vertices in the Cassification Table \ref{tab:table1}. If you start at the vertex $(p+q,p-q)$ in Table \ref{tab:table1}, representing $\G_{p,q}$, and take one step to the left and down, followed by one step to the right and down, you end at the vertex $(p+q+2,p)$ representing $\G_{p+1,q+1}\widetilde{=}M_2(\G_{p,q} )$.   
 Of course in Table \ref{tab:table2}, starting at the Clifford time for $\G_{p,q}$, moving $1$ Clifford hour counterclockwise, followed by $1$ Clifford hour clockwise, brings one back to the same Clifford hour, but the total number of Clifford hours elapsed increases to $n=p+q+2$. 
 
  \begin{thm2} For $k\ge 1$, the geometric algebra $\G_{p+k,q+k}\widetilde{=}M_{2^k}(\G_{p,q} ) $. 
  		\end{thm2}
  		
  	{\bf Proof:} The proof is by induction on $k$. For $k=1$, by definition 
  	\[\G_{p,q}=\R(\be_1, \ldots, \be_p, \bff_1, \ldots, \bff_q) \  {\rm and} \   
  	\G_{p+1,q+1}=\R(\be_1, \ldots, \be_p,\be,  \bff_1, \ldots, \bff_q, \bff ). \]
  	Any element $G\in \G_{p+1,q+1}$ can be expressed in the form
  	\beq  G = g_0 +  g_1\be + g_2 \bff  + g_3 \be \bff,  \label{Gspecbasiskk} \eeq
  	where $g_\mu \in \G_{p,q}$ for $0 \le \mu \le 3$. Applying the matrix standard basis (\ref{mstdbasisG11}) to $\G_{p+1,q+1}$, and noting that 
  	\[  \pmatrix{1 & \be}u_+ \pmatrix{1 \cr \be} =u_+ + \be u_+ \be = u_+ + u_- = 1, \]
  	we calculate
  	
  	\[  G =   \pmatrix{1 & \be}u_+ \pmatrix{1 \cr \be }G \pmatrix{1 & \be}u_+ \pmatrix{1 \cr \be }\] 
  	
  	\[  = \pmatrix{1 & \be}u_+ \pmatrix{G  & G \be \cr  \be  G & \be G \be}u_+  \pmatrix{1 \cr \be}   \]
  	\[ = \pmatrix{1 & \be}u_+ \pmatrix{g_0 + g_3  & g_1-g_2 \cr g_1^- +g_2^-  & g_0^- -g_3^- } \pmatrix{1 \cr \be}=  \pmatrix{1 & \be}u_+[G] \pmatrix{1 \cr \be},  \]
  	where 
  	\[ [G]:=  \pmatrix{g_0 + g_3  & g_1-g_2 \cr g_1^- +g_2^-  & g_0^- -g_3^- } \in M_2(\G_{p,q}), \]
  	and $g^-\equiv\be g \be$ is the operation of inversion in $\G_{p,q}$ obtained by replacing all vectors in $g$ by their negatives. 
  	
  	 Suppose now that the theorem is true for $k-1>1$, so that \[\G_{p+k-1,q+k-1}\widetilde{=}M_{2^{k-1}}(\G_{p,q} ).\] Then for $k$,
  	 \[ \G_{p+k,q+k}\widetilde{=}M_2\big(M_{2^{k-1}}(\G_{p,q})\big)\widetilde{=} M_{2^k}(\G_{p,q} ) ,   \]
  	 which completes the proof. 
  	
  	$\hfill \square$ 
  	
  \no For $p=0=q$, Theorem 2 tells us that
  	$\G_{1,1}\widetilde{=}M_2(\R )$, and for $p=q>1$  
  	\[ \G_{p,p}\widetilde{=}M_2(\G_{p-1,p-1} )\widetilde{=}M_{2^{p}}(\R ).  \]

  	Recalling that $\,^2Q :=\G_{0,3} $, applying {\bf Theorem 1 a)} with $p=0$ and $q=2$ shows that 
  	\beq Q=\G_{0,2} \widetilde{=}\G_{0,3}^+ =\R(\bff_{12}.\bff_{13},\bff_{23}), \label{morequats} \eeq
  	which gives yet another definition of Hamilton's famous quaternions.
  	 To find $\G_{0,4}= \R(\bff_1,\bff_2,\bff_3, \bff_4)$, let $I:=\bff_{123}$, and note that $\bff_4$ and $I$ are 
  	anti-commutative.
    Also, for $\G_{1,1}=\R(\bff_4, \bff_{123})$, and $ \Q :=span_\R\{\bff_{12},\bff_{13},\bff_{23}\}$, then
    $\Q \G_{1,1}=\G_{1,1}\Q$, from which     
  	it follows that 
  	\beq \G_{0,4}{=}\G_{1,1}(\Q)\widetilde{=}M_2(\Q).\label{g04ident} \eeq

The Third Structure Theorem for geometric algebras tells us that in Table \ref{tab:table1} if you start at a vertex representing $\G_{p,q}$, and move four vertices to the left and down, you arrive at the vertex representing $\G_{p+4,q}$, which is isomorphic to
the geometric algebra $\G_{p,q+4}$ four steps to the right and down from $\G_{p,q}$.
On the Clifford Clock in Table \ref{tab:table2}, the relationship follows from the fact that starting at any Clifford time, advancing 4 steps counterclockwise, or clockwise, will bring you to the same Clifford time. 

\begin{thm3} $ \G_{p+4,q}\widetilde=\G_{p,q+4}. \label{bott3} $  

\end{thm3}

{\bf Proof:} Note that
\[ \G_{p+4,q}= \R(\be_a,\be_b,\be_c,\be_d,\be_1, \ldots , \be_p, \bff_1, \ldots, \bff_q  ),    \]
and
\[ \G_{p,q+4}= \R(\be_1, \ldots , \be_p,\bff^\prime_a,\bff^\prime_b,\bff^\prime_c,\bff^\prime_d, \bff_1, \ldots, \bff_q  ),    \]
where 
\[ \bff^\prime_a:= \be_b \be_c \be_d,\ \bff^\prime_b :=  \be_c \be_d \be_a ,\ \bff^\prime_c:= \be_d \be_a \be_b,\ \bff^\prime_d := \be_a \be_b \be_c .  \]
For $s\in \{a,b,c,d\}$, the $\bff^\prime_{s}\in \G_{p+4,q}$ are {\it anticommuting trivectors} which also anticommute with the vector generators of 
$\G_{p,q}$. They also serve as anticommuting vector generators of $ \G_{p,q+4}$, the product of any distinct three of them producing $\pm$ a basis vector in $\G_{p+4,q}$.
 For example, 
 \[ \bff^\prime_a\bff^\prime_b\bff^\prime_c =-\be_ d \in \G_{p+4,q}^1.  \]
 
$\hfill \square $ 

We can now easily verify the validity of the Classification Tables \ref{tab:table1} and \ref{tab:table2}. Our strategy is to first verify the validity of Table \ref{tab:table1} at the vertexes at the points $(p+q,p-q)=(p,-p)$ and $(p+q,p-q)=(p,p)$ corresponding to the
geometric algebras $\G_{0,q}$ and $\G_{p,0}$ for $0 \le p,q \le 7$ along the right and left sides of the triangular table. The validity of the Table at the geometric algebras of interior vertexes then directly follows as a consequence of Theorem 2.     

 We have seen in (\ref{matrixg11}) and (\ref{cmatrixg11}) how to construct the coordinate matrices for the geometric algebras $\G_{1,1}$ and $\G_{1,2}$, respectively.     	
 We have already established that $\G_{0,0}=\R $, $\G_{0,1}=\C $, $\G_{0,2}=\Q $,  and (\ref{g04ident}) establishes that $\G_{0,4}\widetilde{=}M_2(\Q)$.  
 For $\G_{0,5}$, let $i:=\bff_{12345} \in {\cal Z}(\G_{0,5} ) $. Then $i^2=-1$, and from (\ref{cmatrixg11}) it immediately follows that 
 \[ \G_{0,5} =\G_{0,4}(i) \widetilde{=}M_2(\Q)(i)=M_2(\Q(i) )
  =M_4(\C).\] 
 Using Theorem 3, for $\G_{0,6}$,  
 \[ \G_{0,6}\widetilde =\G_{4,2} =\R(\be_1,\be_2,\be_3,\be_4,\bff_1,\bff_2)=  \R(\be_1,\be_{12},\be_{13},\be_{14},\be_1 \bff_1,\be_1 \bff_2)   \]
 \[ \widetilde{=}\G_{3,3} \widetilde{=}M_8(\R).  \]
 In the last two steps, we have used both Theorem 1 b) and Theorem 2, respectively. Finally, for $\G_{0,7}$, note that $I:=\bff_{1\cdots 7}$ has the property that
 $I \in {\cal Z}(\G_{0,7} )$ and $I^2=1$. It follows from (\ref{doubleG}) that 
 \[  \G_{0,7}=\G_{0,6}(I)=\,^2\G_{0,6} \widetilde{=}\, ^2M_8(\R)=M_8(\,^2\R).  \]  
 
 We now turn our attention to $\G_{p,0}$ for $0\le p \le 7$. We have already established that $\G_{1,0}\widetilde{=}\,^2\R $. Using Theorem 1 b), and Theorem 2, 
 \[ \G_{2,0}\widetilde=\G_{1,1}\widetilde{=}M_2(\R). \]
 For $\G_{3,0}$, using Theorem 1 b),
 \[  \G_{3,0} \widetilde{=}\G_{1,2}\widetilde{=}M_2(\C).   \]  
 By Theorem 1 b),
 \[ \G_{4,0} \widetilde{=}\G_{1,3}\widetilde{=}\G_{1,2}(\bff_3)\widetilde{=}M_2(\C)(\bff_3)\widetilde{=}M_2(\Q). \] 
 Similarly,
   \[ \G_{5,0} \widetilde{=}\G_{1,4}\widetilde{=}\G_{1,3}(\bff_4)\widetilde{=}M_2(\C)(\bff_4)\widetilde{=}M_2(\,^2\Q), \]
   \[ \G_{6,0} \widetilde{=}\G_{1,5}\widetilde{=}\G_{1,4}(\bff_5)\widetilde{=}M_2(\,^2\Q)(\bff_5)\widetilde{=}M_4(\Q), \]  
   	and, finally,
    \[ \G_{7,0} \widetilde{=}\G_{1,6}\widetilde{=}\G_{1,5}(\bff_6)\widetilde{=}M_4(\Q)(\bff_6)\widetilde{=}M_8(\C). \]

  \section{Geometric algebras $\G_{k,k}$ and $\G_{k,k+1}$   }
  Coordinate geometric matrices for the geometric algebras $\G_{1,1}$ and $\G_{1,2}$ were constructed in {\bf Section 3}. We now show how the construction given in
  (\ref{mstdbasisG11}) and (\ref{cmatrixg11}) generalizes to $\G_{k,k}$ and $\G_{k,k+1}$ for $k>1$. For both $\G_{1,1}$ and $\G_{1,2}$, the standard vector and null vector basis is given by  
  \[ \pmatrix{1 \cr \be_1 }u_1\pmatrix{1 & \be_1}   =\pmatrix{1 \cr \ba_1 }u_1\pmatrix{1 & \bb_1},  \]
  respectively, where $u_1=\frac{1}{2}(1+\be_1\bff_1)=\bb_1\ba_1$. Given $g\in \G_{1,1}$ or
  $g\in \G_{1,2}$, respectively, for the real or complex matrix $[g]$ of $g$, respectively, 
  \[ g = \pmatrix{1 & \be_1 }u_1[g]\pmatrix{1 \cr \be_1}   =\pmatrix{1 & \ba_1 }u_1[g]\pmatrix{1 \cr \bb_1}. \]
  
  For both $\G_{k,k}$ and $\G_{k,k+1}$, utilizing the {\it directed Kronecker product} \cite[p.82]{matrixG}, the standard and null basis are given by  
  \beq \pmatrix{1 \cr \be_1 }\overrightarrow{\otimes}\cdots\overrightarrow{\otimes}\pmatrix{1 \cr \be_k }u_{1\cdots k}\pmatrix{1 & \be_k} \overleftarrow{\otimes}\cdots\overleftarrow{\otimes}\pmatrix{1 & \be_1 } \label{kstdbasis} \eeq 
  \beq =\pmatrix{1 \cr \ba_1 }\overrightarrow{\otimes}\cdots\overrightarrow{\otimes}\pmatrix{1 \cr \ba_k }u_{1\cdots k}\pmatrix{1 & \bb_k} \overleftarrow{\otimes}\cdots\overleftarrow{\otimes}\pmatrix{1 & \bb_1 },  \label{knullbasis} \eeq
  respectively, where $u_j=\frac{1}{2}(1+\be_j\bff_j)=\bb_j\ba_j$, $u_{1\cdots k}=
  u_1\cdots u_k$, and, analogous to (\ref{nullabdef}),  
  \beq \ba_j := \frac{1}{2}(\be_j+\bff_j), 
  \ {\rm and} \ \bb_j :=\frac{1}{2}(\be_j -\bff_j). \label{nullabjdef} \eeq
  
  We now illustrate the use of the null standard basis (\ref{nullstdbasis}) for $k=2$.  Given $g\in \G_{2,2}$ or
  $g\in \G_{2,3}$, respectively, for the real or complex matrix $[g]$ of $g$, respectively, 
  \[ g = \pmatrix{1 & \ba_1 }\overrightarrow{\otimes}\pmatrix{1&\ba_2 } u_{12}[g]\pmatrix{1 \cr \bb_2}\overleftarrow{\otimes}\pmatrix{1 \cr \bb_1 }  \] 
  \[ =\pmatrix{1 & \ba_1 & \ba_2 & \ba_{12} }u_{12} [g] \pmatrix{1 \cr \bb_1 \cr \bb_2 \cr \bb_{21}   }  , \]
  with respect to the null standard basis for $k=2$, given by
  \[ \pmatrix{1 \cr \ba_1 }\overrightarrow{\otimes}\pmatrix{1 \cr \ba_2 }u_{12}\pmatrix{1 & \bb_2} \overleftarrow{\otimes}\pmatrix{1 & \bb_1 } =\pmatrix{1 \cr \ba_1 \cr \ba_2 \cr \ba_{12} }u_{12}  \pmatrix{1 & \bb_1 & \bb_2 & \bb_{21}   } \] 
   \[ =  \pmatrix{u_{12} & \bb_{1} u_2 & \bb_{2}u_1 & \bb_{21} 
   	\cr \ba_1 u_2 & u_1^\dagger u_2 & \ba_1 \bb_2 &-\bb_2 u_1^\dagger   \cr
   	\ba_{2}u_1  & \ba_2 \bb_1 & u_1 u_2^\dagger & \bb_1 u_2^\dagger  \cr
   	\ba_{12} & -\ba_2 u_1^\dagger & \ba_1u_2^\dagger & u_{12}^\dagger}.   \]  
   Further details of the construction are given in \cite[pp\,82-88]{matrixG}, and tables are given for the null standard basis of the geometric algebras $\G_{k,k}$ and $\G_{k,k+1}$ up to $k=5$. 
   
   The classification of the complex geometric algebras  
   $ \G_{n}(\C)=\G_{p,q}(\C)$, for $n=p+q$, greatly simplifies because
   \[ \G_{p,q}(\C)=\R(\be_1,\ldots,\be_p,\bff_1,\ldots,\bff_q)(\C)= \C(\be_1,\ldots,\be_p,\bff_1,\ldots,\bff_q)        \]  
   \[ =\C(\be_1,\ldots,\be_p,i\bff_1,\ldots,i\bff_q)=\G_n(\C).     \]  
 The Classification Table for complex geometric algebra is given in Table 4, \cite[p\,75]{matrixG}.

  \begin{table}
  	\centering
  	
  	\caption{Classification of Complex Geometric Algebras.}
  	
  	\centering
  	\label{tab:table4}
  	\begin{tabular}{cccccccc} \\
  		$\G_{0}(\C)$ &$\G_{1}(\C)$  &$\G_{2}(\C)$   &$\G_{3}(\C)$ 
  		&$\G_{4}(\C)$    &$\G_{5}(\C)$  & $\G_{6}(\C)$  &$\G_{7}(\C)$    \\ 
  		$\C $  & $^2\C $& $M_2(\C)$ & $ M_2(^2\C ) $ & $M_4 (\C )$ & $M_4(^2\C)$ & $M_8(\C ) $& $M_8(^2\C) $ 
  	\end{tabular}  
  \end{table}

	\end{document}